\documentclass[11pt]{amsart}
\usepackage{amssymb, latexsym}
\theoremstyle{plain}
\newtheorem{theorem}{Theorem}
\newtheorem{corollary}{Corollary}

\newtheorem {lemma}{Lemma}
\newtheorem{proposition}{Proposition}

\theoremstyle{remark}

\numberwithin{equation}{section}

\begin{document}

\title[]%
 {Transience/Recurrence and the speed of  a one-dimensional random walk  in a ``have your cookie and eat it'' environment}

\author{Ross G. Pinsky}
\address{Department of Mathematics\\
Technion---Israel Institute of Technology\\
Haifa, 32000\\ Israel} \email{pinsky@math.technion.ac.il}
\urladdr{http://www.math.technion.ac.il/~pinsky/}

\date{}

\begin{abstract}
Consider a variant of the simple  random walk on the integers, with the following transition mechanism. At each site $x$,  the probability of jumping
to the right is $\omega(x)\in[\frac12,1)$,  until the first time the process jumps to the left from site $x$, from which time onward
the probability of jumping to the right is $\frac12$. We investigate the transience/recurrence properties of this process
in both deterministic and stationary, ergodic environments $\{\omega(x)\}_{x\in Z}$.
In deterministic environments, we also study the speed of the process.
\end{abstract}
\maketitle
\section{Introduction and Statement of Results}
A simple random walk on the integers in a cookie environment \newline $\{\omega(x,i)\}_{x\in Z, i\in\mathcal{N}}$ is a
stochastic process $\{X_n\}_{n=0}^\infty$ defined on the integers  with the following transition mechanism:
$$
P(X_{n+1}=X_n+1|\mathcal{F}_n)=1-P(X_{n+1}=X_n-1|\mathcal{F}_n)=\omega(X_n,V(X_n;n)),
$$
 where
$\mathcal{F}_n=\sigma(X_0,\cdots,X_n)$ is the filtration of the process up to time $n$
 and $V(x;n)=\sum_{k=0}^n\chi_x(X_k)$.
A graphic description can be given as follows. Place an infinite sequence of cookies at each site $x$.
The $i$-th time the process reaches the site $x$, the process eats the $i$-th cookie
at site $x$, thereby empowering it to  jump to the right
with probability $\omega(x,i)$ and to the left with probability $1-\omega(x,i)$.
If $\omega(x,i)=\frac12$, the corresponding cookie is a placebo; that is, its effect on the direction of the random walk
is neutral. Thus, when for each $x\in Z$ one
has $\omega(x,i)=\frac12$, for all sufficiently large $i$, one can consider  the process  as a simple, symmetric random walk with
 a finite number of symmetry-breaking cookies per site.

Letting $P_x$ denote probabilities for the process conditioned on $X_0=x$, and letting
Define the process
to be  transient or recurrent respectively according to whether  $P_y(X_n=x\ \text{i.o.})$
is equal to 0 or 1 for all $x$ and all $y$. Of course from this definition, it is not a priori clear that the process is either
transient or recurrent.
When there is exactly one cookie per site and
$\omega(x,1)=p\in(\frac12,1)$, for all $x\in Z$,  the random walk has been called an excited random walk.
It is not hard to show that this random walk is recurrent  \cite{D}, \cite{BW}.
Zerner showed \cite{Z} that when multiple cookies are allowed per site, then it is possible to have transience. For
 a   stationary, ergodic  cookie environment
$\{\omega(x,i)\}_{x\in Z, i\in\mathcal{N}}$, where  the cookie values are in $[\frac12, 1)$,
he gave  necessary and sufficient conditions for transience or recurrence.
In particular, in the deterministic case where
there are $k\ge2$ cookies per site and $\omega(x,i)=p$, for $i=1,\cdots, k$ and all $x\in Z$,  a phase transition
with regard to transience/recurrence occurs at  $p=\frac12+\frac1{2k}$: the process is recurrent
if $p\in(\frac12,\frac12+\frac1{2k}]$ and transient if $p>\frac12+\frac1{2k}$.
Zerner also showed in the above ergodic setting
that if there is no excitement after the second visit---that is, with probability one, $\omega(0,i)=\frac12$, for $i\ge3$---then
the random walk is non-ballistic; that is, it has 0 speed: $\lim_{n\to\infty}\frac{X_n}n=0$ a.s..
He left as an open problem the question of whether it is possible to have positive speed when there are
a bounded number of cookies per site.
Recently, Basdevant and Singh  \cite{BS} have given an explicit criterion for determining whether the speed is positive or zero
in the case of a deterministic, spatially homogeneous environment with a finite number of cookies.
In particular, if there are $k\ge3$ cookies per site and $\omega(x,i)=p$, for $i=1,\cdots, k$ and all $x\in Z$, then a phase transition with regard to ballisticity
occurs at $p=\frac12+\frac1k$:
one has $\lim_{n\to\infty}\frac{X_n}n=0$ a.s. if $p\in(\frac12,\frac12+\frac1k]$ and
 $\lim_{n\to\infty}\frac{X_n}n>0$ a.s. if $p>\frac12+\frac1k$. In particular, this shows
 that it is possible to have  positive speed with three cookies per site.
 Kosygina and Zerner \cite{KZ} have considered cookie environments that can take values in $(0,1)$, so that
the cookie bias can be to the right at certain sites and to the left at other sites.

In this paper we consider one-sided cookie bias in  what we dub a ``have your cookie and eat it environment.'' There
is one cookie at each site $x$ and its value is $\omega(x)\in[\frac12,1)$. If the random walk moves to the right from site $x$,
then the cookie at site $x$ remains in place unconsumed, so that the next time the random walk reaches  site $x$, the cookie
environment is still in effect. The cookie at site $x$  only gets consumed when the random walk moves leftward from $x$.
After this occurs, the transition mechanism at  site $x$ becomes that of   the  simple, symmetric random walk.
We study the transience/recurrence and ballisticity of such   processes.
We have several motivations for studying this process. One of the motivations for studying the cookie random walks
described above is to understand processes with self-interactions. The random walk in the
``have your cookie and eat it'' environment adds an extra level of self-interaction---in previous works, the number of times
cookies are employed at a given site is an external parameter, whereas here it depends on the behavior of the process.
Another novelty is that we are able to calculate certain quantities
explicitly, quantities that have not been calculated explicitly in  previous models.
For certain sub-ranges of the appropriate parameter, we
 calculate  explicitly the speed   in the ballistic case
and the probability of ever reaching 0 from 1  in the transient case. Finally, although  we use some of the ideas and methods
of \cite{Z} and \cite{BS} at the beginning of our proofs, we  introduce some new techniques, which may well
 be useful  in other cookie problems also.

When considering a deterministic environment $\omega$, we will denote probabilities and expectations  for the process starting at $x$ by $P_x$
and $E_x$ respectively, suppressing the dependence on the environment, except in Proposition \ref{mono} below, where we
need to compare probabilities for two different environments.
When considering
a stationary, ergodic environment  $\{\omega(x)\}_{x\in Z}$,
we will denote the probability measure on the environment by $\mathcal{P}^S$ and the corresponding
expectation by $\mathcal{E}^S$.
Probabilities for the process  in a fixed environment $\omega$
 starting from $x$ will be denoted
by $P^\omega_x$ (the \it quenched\rm\ measure).

Before stating the results, we present a lemma concerning transience and recurrence.
Let
\begin{equation}
T_x=\inf\{n\ge0:X_n=x\}.
\end{equation}

\begin{lemma}\label{lemma}
Let $\{X_n\}_{n=0}^\infty$ be a random walk in a ``have your cookie and eat it'' environment $\{\omega(x)\}_{x\in Z}$,
where $\omega(x)\in[\frac12,1)$, for all $x\in Z$.

\noindent 1. Assume that the environment $\omega$ is deterministic and periodic: for some $N\ge1$, $\omega(x+N)=\omega(x)$, for all $x\in Z$.

i. If  $P_1(T_0=\infty)=0$,
 then
the process is recurrent.

ii. If
$P_1(T_0=\infty)>0$, then the process is transient and $\lim_{n\to\infty}X_n=\infty$ a.s..

\noindent 2. Assume that  the environment $\omega$ is
stationary and ergodic.
Then either $P^\omega_1(T_0=\infty)=0$ for $\mathcal{P}^S$-almost every $\omega$ or
$P^\omega_1(T_0=\infty)>0$ for $\mathcal{P}^S$-almost every $\omega$.

 i. If  $P^\omega_1(T_0=\infty)=0$ for $\mathcal{P}^S$-almost every $\omega$,
 then for $\mathcal{P}^S$-almost every $\omega$
the process is  recurrent.

 ii. If $P^\omega_1(T_0=\infty)>0$ for $\mathcal{P}^S$-almost every $\omega$,
then for $\mathcal{P}^S$-almost every $\omega$ the process is  transient and
$\lim_{n\to\infty}X_n=\infty$ a.s..
\end{lemma}

The following theorem concerns transience/recurrence in deterministic environments.

\begin{theorem}\label{t/r-det}
Let $\{X_n\}_{n=0}^\infty$ be a random walk in a deterministic ``have your cookie and eat it'' environment $\{\omega(x)\}_{x\in Z}$,
where $\omega(x)\in[\frac12,1)$ for all $x\in Z$.
\begin{enumerate}
\item[i.]
Let $\omega(x)=p\in(\frac12,1)$, for all $x\in Z$.
Then
$$
\begin{cases}
P_1(T_0=\infty)=
0,\ \text{if}\ p\le \frac23;\\
\frac{3p-2}p\le P_1(T_0=\infty)\le \frac1p\frac{3p-2}{2p-1},\ \text{if}\ p\in(\frac23,1).
\end{cases}
$$
In particular the process is recurrent if $p\le\frac23$ and transient if $p>\frac23$.
\item[ii.] Let $\omega$ be periodic with period $N>1$.
Then the process is recurrent if
$\frac1N\sum_{x=1}^N\frac{\omega(x)}{1-\omega(x)}\le2$ and transient if
$\frac1N\sum_{x=1}^N\frac{\omega(x)}{1-\omega(x)}>2$.

\end{enumerate}
\end{theorem}
\bf\noindent Remark 1.\rm\ We note the following  heuristic connection between part (i) of the theorem and Zerner's result
that was described above. Zerner proved that when there are $k$ cookies per site, all with strength $p\in(\frac12,1)$, then
the process is recurrent if  $p\le\frac12+\frac1{2k}$ and transient if $p>\frac12+\frac1{2k}$. For the process
in a ``have your cookie and eat it'' environment, if it is recurrent so that it continually returns to every site, then
for each site, the expected number of times it will jump from that site according to the probabilities $p$ and $1-p$ is
$\sum_{m=1}^\infty mp^{m-1}(1-p)=\frac1{1-p}$.
Therefore, on the average, the process in  a ``have your cookie and eat it'' environment with parameter $p$  behaves like
a $k$ cookies per site process with $k=\frac1{1-p}$ (which of course is usually not an integer).
If one substitutes this value of $k$ in the equation $p=\frac12+\frac1{2k}$ and solves for $p$,
one obtains $p=\frac23$.
\medskip

\noindent \bf Remark 2.\rm\ In the case of a periodic environment, it follows from Jensen's inequality
that if the average cookie strength  $\frac1N\sum_{x=1}^N\omega(x)>\frac23$, then \newline
$\frac1N\sum_{x=1}^N\frac{\omega(x)}{1-\omega(x)}>2$, and thus the process is transient. Indeed, it is possible
to have transience with the average cookie strength arbitrarily close to $\frac12$: let $\omega(x)=\frac12$, for $x=1,\cdots, N-1$, and
$\omega(N)>\frac{N+1}{N+2}$.

\medskip
\noindent \bf Remark 3.\rm\
Part (i) of  Theorem \ref{t/r-det}  gives two-sided bounds on $P_1(T_0=\infty)$.
In Proposition \ref{T0} in section 4, we will prove that
$P_1(T_0=\infty)\ge\frac{3p-2}{2p-1}$, for $p>\frac23$,  with equality for $p\ge\frac34$.
We believe that equality holds for all $p>\frac23$.

\medskip

The next
theorem treats transience/recurrence in stationary, ergodic environments  $\{\omega(x)\}_{x\in Z}$.
\begin{theorem}\label{t/r-random}
Let $\{X_n\}_{n=0}^\infty$ be a random walk in a stationary, ergodic ``have your cookie and eat it'' environment $\{\omega(x)\}_{x\in Z}$.
Then for $\mathcal{P}^S$-almost every environment $\omega$,  the process is  $P^\omega_x$-recurrent if  $\mathcal{E}^S\frac{\omega(0)}{1-\omega(0)}\le2$ and
 $P^\omega_x$-transient if $\mathcal{E}^S\frac{\omega(0)}{1-\omega(0)}>2$.
\end{theorem}

We now turn to the speed of the random walk.
\begin{theorem}\label{ball}
Let $\{X_n\}_{n=0}^\infty$ be a random walk in a deterministic ``have your cookie and eat it'' environment
with $\omega(x)=p$, for all $x$.

\noindent
If $p<\frac34$, then the process is non-ballistic; one has
 $\lim_{n\to\infty}\frac{X_n}n=0$ a.s..

\noindent If $p>\frac34$, then the process is ballistic; one has
$\lim_{n\to\infty}\frac{X_n}n>0$ a.s..
\newline In fact,
\begin{equation}\label{explicit}
\lim_{n\to\infty}\frac{X_n}n\ge 4p-3,\  \text{with equality if}\ \ \frac{10}{11}\le p<1.
\end{equation}

\end{theorem}
\bf\noindent Remark 1.\rm\ We believe that equality holds in \eqref{explicit} for all $p\in(\frac34,1)$.
In particular, this would prove non-ballisticity in the borderline case $p=\frac34$.
\medskip

\bf\noindent Remark 2.\rm\ We note the following heuristic connection between  Theorem \ref{ball} and
the result of Basdevant and Singh that was described above. That result  states that
if there are
$k$ cookies per site, all with strength $p\in(\frac12,1)$, then
the process has 0 speed if  $p\le\frac12+\frac1{k}$ and positive speed if $p>\frac12+\frac1{k}$.
As noted in the remark following Theorem \ref{t/r-det}, on the average the process in a ``have your cookie and eat it'' environment
behaves like a $k$ cookie per site process with $k=\frac1{1-p}$. If one substitutes this value of $k$ in the equation
$p=\frac12+\frac1k$ and solves for $p$, one obtains $p=\frac34$.

\medskip

In section 2 we prove a monotonicity result with respect to environments and the use it to prove Lemma \ref{lemma}.
The proofs of  Theorems \ref{t/r-det} and \ref{t/r-random} are given in section 3, and the proof
of Theorem \ref{ball} is given  in  section 4.

\section{A monotonicity Result and the Proof of Lemma \ref{lemma}}
A standard coupling \cite[Lemma 1]{Z} shows that
$\{X_n\}_{n=0}^\infty$ under $P^{\omega}_y$ stochastically dominates the simple, symmetric random walk starting from $y$.
However,   standard couplings don't work for comparing  random walks in two different
``have your cookie and eat it'' environments (or multiple cookie environments)---see \cite{Z}.
Whereas standard couplings are ``space-time'' couplings, we will construct a ``space only'' coupling, which allows the time
parameter of the two processes to differ.

\begin{proposition}\label{mono}
Let $\omega_1$ and $\omega_2$ be environments with $\omega_1(z)\le \omega_2(z)$, for all $z\in Z$. Denote probabilities
for the random walk in the  ``have your cookie and eat it'' environment  $\omega_i$ by $P_{\cdot}^{\omega_i}$, $i=1,2$.
Then
\begin{equation}\label{comp}
P^{\omega_1}_y(T_z\le T_x\wedge N)\le P^{\omega_2}(T_z\le T_x\wedge N), \ \text{for} \ x<y<z\ \text{and}\ N>0.
\end{equation}
In particular then, $P^{\omega_1}_y(T_z< T_x)\le P^{\omega_2}(T_z< T_x)$ and
$P_y^{\omega_1}(T_x=\infty)\le P_y^{\omega_2}(T_x=\infty)$.
\end{proposition}
\bf\noindent Remark. \rm\ The same result was proven in \cite{Z} for one-side multiple cookies,
where one assumes that the two environments satisfy $\omega_1(x,i)\le \omega_2(x,i)$,
for all $x\in Z$ and all $i\ge1$. The proof does not  use coupling.
Our proof extends to that setting if one assumes that the environments satisfy
$\omega_1(x,j)\le\omega_2(x,i)$, for all $x\in Z$ and all $j\ge i\ge1$.

\begin{proof}
On a probability space $(\Omega, \mathcal{P})$, define a sequence $\{U_k\}_{k=1}^\infty$  of IID  random variables uniformly
distributed on $[0,1]$.
Let $\mathcal{F}_k=\sigma(U_j:j=1,\cdots, k)$.
We now couple two processes, $\{X_n^1\}$ and $\{X_n^2\}$, on $(\Omega,\mathcal{P})$ so
that $\{X_n^i\}$ has the distribution of the  $P^{\omega_i}_y$ process, $i=1,2$.
The coupling will be for all times $n$, without reference to $T_x$ or $T_z$ or $N$ in the statement of the proposition.
 We begin with $X_0^1=X_0^2=y$.
 For $i=1,2$, we let $X_1^i=y+1$ or $X_1^i=y-1$ respectively depending on whether $U_1\le \omega_i(y)$ or $U_1>\omega_i(y)$.
 If $U_1\le \omega_1(y)$ or $U_1>\omega_2(y)$, then both processes moved together. In this case, use $U_2$ to continue
 the coupling another step. Continue like this until the first time $n_1\ge1$ that the processes differ.
 Necessarily,
 $X_{n_1}^1=X_{n_1}^2-2$. Up until this point the coupling has used the random variables $\{U_k\}_{k=1}^{n_0}$.
 Now leave the process $\{X_n^2\}$ alone and continue moving the process $\{X^1_n\}$, starting with the random variable
 $U_{n_0+1}$. By comparison with the simple, symmetric random walk, the $\{X_n^1\}$ process will eventually return to the level
 $X_{n_0}^2$ of the temporarily frozen $\{X_n^2\}$ process. Denote this time by $m_0$. So $X_{m_0}^1=X_{n_0}^2$.
 Now start moving both of the processes together again, starting with the first unused random variable; namely, $U_{m_0+1}$.
 Continue the coupling until the process differ again. Then as before, run only the $\{X_n^1\}$ process until it again reaches
 the newly, temporarily frozen level of the $\{X_n^2\}$ process, etc.,etc..

 Clearly, the $\{X_n^i\}$ process has the distribution of the $P^{\omega_i}_y$ process, $i=1,2$.
Let $T^i_r$ denote the hitting time of $r$ for $\{X_n^i\}$. We claim that
$\{T^1_z\le T^1_y\wedge N\}\subset\{T^2_z\le T^2_x\wedge N\}$, from which it follows that
$\mathcal{P}(T^1_z\le T^1_x\wedge N)\le \mathcal{P}(T^2_z\le T^2_x\wedge N)$, thereby proving \eqref{comp}.
We explain this, using the notation in the previous paragraph. Consider the processes
$\{X^i_n\}$ restricted to $\mathcal{F}_{n_0}$.  From the construction,
it is clear that
$\{T^1_z\le T^1_y\wedge N\}\cap\{T^1_z\le n_0\}=\{T^2_z\le T^2_x\wedge N\}\cap\{T^2_z\le n_0\}$.
If this event does not occur, then
 consider the processes restricted to $\mathcal{F}_{m_0}$.
Between time $n_0$ and $m_0$ the $\{X^1_n\}$ process might hit
$x$. Because of this, and because it is possible that  $m_0\ge N$, it follows that
$\{T^1_z\le T^1_y\wedge N\}\cap\{T^1_z\le m_0\}\subset\{T^2_z\le T^2_x\wedge N\}\cap\{T^2_z\le m_0\}$.
Continuing like this on each  succeeding pair of intervals, the first one of which  the processes move together in a coupled manner
and the second one of which the $\{X_n^2\}$  process is frozen and the $\{X_n^1\}$ process moves alone until it returns to the level
of the frozen $\{X_n^2\}$ process,
 we conclude that
 $\{T^1_z\le T^1_y\wedge N\}\subset\{T^2_z\le T^2_x\wedge N\}$.

 The last line of the proposition follows by letting $N\to\infty$ and then letting $z\to\infty$.

\end{proof}

\noindent\bf Proof of Lemma \ref{lemma}.\it\ (1-i).\rm\ Since
$P_1(T_0=\infty)=0$, the process starting at 1 will eventually hit 0. It will then return to 1, by comparison with the simple,
symmetric random walk. Upon returning to 1, the  current environment, call it $\omega'$, satisfies $\omega'\le\omega$, since
a jump to the left along the way at any site $x$ changed the jump mechanism at that site from $\omega(x)$ to $\frac12$.
By Proposition \ref{mono}, it then follows that
the process will hit 0  again. This process continues indefinitely; thus,  $P_1(X_n=0\ \text{i.o.})=1$.
For any $x\in Z$, there exists a $\delta_x>0$ such that regardless of what the current environment is, if the process is at 0 it will
go directly to $x$ in $|x|$ steps. Thus from $P_1(X_n=0\ \text{i.o.})=1$ it follows that
$P_1(X_n=x\ \text{i.o.})=1$, for all $x\in Z$. Now let $y>1$ and let $x<y$. There is a positive probability that the process  starting at 1
will
make it first  $y-1$ steps to the right, ending up at $y$ at time $y-1$. The current environment will then still be $\omega$.
Since $P_1(X_n=x\ \text{i.o.})=1$, it then follows that $P_y(X_n=x\ \text{i.o.})=1$.
Now let $y\le 0$ and $x<y$. By comparison with the simple, symmetric random walk, the process starting at $y$ will eventual
hit 1. The current environment now, call it $\omega'$, satisfies $\omega'\le \omega$. Thus
since $P_1(X_n=x\ \text{i.o,})=1$, it follows from Proposition \ref{mono} that $P_y(X=x\ \text{i.o.})=1$.
\ \ \ \

\ \

\noindent \it (1-ii).\rm\ Let $P_1(T_0=\infty)>0$. We first show that
$P_y(T_x=\infty)>0$, for all $x<y$.
To  do this, we   assume that for some $x_0<y_0$ one has $P_{y_0}(T_{x_0}=\infty)=0$, and then
 come to a  contradiction.
By
 the argument in (1-i), it follows that $P_{y_0}(X_n=x_0\ \text{i.o})=1$, and then that
$P_{y_0}(X_n=x\ \text{i.o.})=1$, for all $x\in Z$. In particular, $P_{y_0}(X_n=0\ \text{i.o.})=1$.
If $y_0=1$, we have arrived at a contradiction. Consider now the case $y_0<1$. Since there is a positive probability of going directly from
$y_0$ to 1 in the first $1-y_0$ steps, in which case the current environment will still be $\omega$, and since
$P_{y_0}(X_n=0\ \text{i.o.})=1$, it follows that
$P_1(X_n=0\ \text{i.o.})=1$, a contradiction. Now consider the case $y_0>1$. Since $P_1(T_0=\infty)>0$, there is a positive probability
that starting from 1 the process will hit $y_0$ without having hit 0, and then after reaching $y_0$, continue and  never ever hit 0.
But when this process reached $y_0$, the current environment, call it $\omega'$, satisfied $\omega'\le\omega$.
Since $P_{y_0}(X_n=0\ \text{i.o.})=1$, it follows from Proposition \ref{mono} that after hitting $y_0$ the process would in fact hit 0 with
probability 1, a contradiction.

From the fact that
 $P_y(T_x=\infty)>0$, for all $x<y$, and from the
 periodicity of the environment, it follows that $\gamma\equiv\inf_{x\in Z} P_x(T_{x-1}=\infty)=\inf_{x\in\{0,1,\cdots, N-1\}}P_x(T_{x-1}=\infty)>0$.
From this we can prove transience. Indeed, fix $x$ and $y$. By  comparison with the simple, symmetric random walk, with $P_y$-probability
one, the process
will attain a new maximum  infinitely often. The current environment at and to the right of any new  maximum
coincides with $\omega$. Thus,
 with probability at least $\gamma$
the process  will never again  go below that maximum. From this it follows that $P_y(X_n=x\ \text{i.o})=0$ and that
 $\lim_{n\to\infty}X_n=\infty$ a.s..
\ \ \

\ \

\noindent\it (2).\rm\ It follows from the proof of  part (1) that if $P^\omega_1(T_0=\infty)=0$, then $P^\omega_{x+1}(T_x=\infty)=0$ for all $x\ge1$, and if
$P^\omega_1(T_0=\infty)>0$, then $P^\omega_{x+1}(T_x=\infty)>0$ for all $x\ge1$.
From this and
 the  stationarity and ergodicity assumption, it follows  that either $P^\omega_1(T_0=\infty)=0$ for $\mathcal{P}^S$-almost every $\omega$
or $P^\omega_1(T_0=\infty)>0$ for $\mathcal{P}^S$-almost every $\omega$.
For an $\omega$ such that
 $P^\omega_1(T_0=\infty)=0$, the
  proof of recurrence is the same as the corresponding  proof in 1-i.
For an $\omega$ such that $P^\omega_1(T_0=\infty)>0$, the proof of
transience and the proof
that $\lim_{n\to\infty}X_n=\infty$ a.s.
are similar to the corresponding proofs in 1-ii.  In 1-ii
we used the fact that $\inf_{x\in Z}P_{x+1}(T_x=\infty)>0$.
However, it is easy to see that all we needed there is
$\limsup_{x\to\infty}P_{x+1}(T_x=\infty)>0$. By the stationarity
and ergodicity, this condition holds for a.e. $\omega$.

\hfill$\square$


\section{Proofs of Theorems \ref{t/r-det} and \ref{t/r-random} }
\noindent \bf Proof of Theorem \ref{t/r-det}.\rm\
For $x\in Z$, define $\sigma_x=\inf\{n\ge1: X_{n-1}=x, X_n=x-1\}$. Let $D_0^x=0$ and for $n\ge1$, let
\begin{equation}
D_n^x=\#\{m<n:  X_m=x, \sigma_x> m\}.
\end{equation}
Define
\begin{equation}
D_n=\sum_{x\in Z} (2\omega(x)-1)D_n^x.
\end{equation}
It follows easily from the definition of the random walk in a ``have your cookie and eat in'' environment that
$X_n-D_n$ is a martingale. Doob's optional stopping theorem gives
\begin{equation}\label{Doob}
nP_1(T_n<T_0)=1+E_1D_{T_0\wedge T_n}=1+\sum_{x=1}^{n-1}(2\omega(x)-1)E_1D^x_{T_0\wedge T_n}.
\end{equation}
One has $\lim_{n\to\infty}P_1(T_n<T_0)=P_1(T_0=\infty)$.
 Let $\gamma_x=P_{x+1}(T_x=\infty)$.

Assume for now  that the process is transient.
For the meantime, we do not distinguish between parts (i) and (ii) of the theorem; that is, between $N=1$ and $N>1$.
By the transience and   Lemma \ref{lemma},
  $\gamma_x>0$ for all $x$.
Consider the random variable $D_{T_0}^x$ under $P_1$, for $x\ge1$. We now construct a random variable which
dominates $D_{T_0}^x$. If $T_0<T_x$, then $D_n^x=0$.
If $T_x<T_0$, and at time $T_x$ the process jumps to $x-1$, which occurs  with probability
$1-\omega(x)$, or jumps to $x+1$ and then never returns to
$x$, which occurs with probability $\omega(x)\gamma_x$, then one has $D_{T_0}^x=1$. Otherwise, $D_{T_0}^x>1$.
In this latter case, the process will return to $x$ from the right.
Upon returning, the process jumps to $x-1$ with probability $1-\omega(x)$, in which case $D_{T_0}^x=2$.
Whenever the process jumps again to $x+1$, it may not ever return to $x$, but
for the domination, we ignore this possibility and  assume that it returns  to $x$ with probability 1.
Taking the above into consideration, it follows that
under $P_1$, the random variable $D_{T_0}^x$ is stochastically dominated by
$I_xZ_x$, where $I_x$ and $Z_x$ are independent random variables satisfying
\begin{equation}\label{IZ}
\begin{aligned}
&P(I_x=1)=1-P(I_x=0)=P_1(T_x<T_0),\\
&
\begin{cases} P(Z_x=1)=1-\omega(x)+\omega(x)\gamma_x;\\
 P(Z_x=k)=(\omega(x)-\omega(x)\gamma_x)\omega^{k-2}(x)(1-\omega(x)),\ k\ge2.
\end{cases}
\end{aligned}
\end{equation}
Thus,  $E_1D_{T_0}^x\le EI_xZ_x=P_1(T_x<T_0)EZ_x$. One has  $EZ_x=\frac{1-\omega(x)\gamma_x}{1-\omega(x)}$. Thus,
\begin{equation}\label{Dx}
E_1D^x_{T_0}\le \frac{1-\omega(x)\gamma_x}{1-\omega(x)}P_1(T_x<T_0).
\end{equation}
Substituting \eqref{Dx} into \eqref{Doob} and using the monotonicity of $D_n$ gives
\begin{equation}\label{key}
P_1(T_n<T_0)\le\frac1n+\frac1n\sum_{x=1}^{n-1}(2\omega(x)-1)\frac{1-\omega(x)\gamma_x}{1-\omega(x)}P_1(T_x<T_0).
\end{equation}

Consider now   part (i) of the theorem. In this case $\omega(x)=p$ for all $x$, and
thus $\gamma\equiv\gamma_x$ is independent of $x$. Letting
$n\to\infty$ in \eqref{key} and using the transience assumption, $\gamma>0$,
one obtains
\begin{equation*}
\frac{(2p-1)(1-p\gamma)}{1-p}\ge1,
\end{equation*}
or equivalently,
\begin{equation}\label{finalupper-i}
\gamma=P_1(T_0=\infty)\le \frac1p\frac{3p-2}{2p-1}.
\end{equation}
Since \eqref{finalupper-i} was derived under the assumption that $\gamma>0$, it follows that $p>\frac23$ is  a necessary condition
for transience.
Note
also that \eqref{finalupper-i} gives the upper bound on $P_1(T_0=\infty)$ in part (i) of the theorem.

Now consider part (ii) of the theorem.
In this case $\omega(x)$ and thus also $\gamma_x$ are periodic with period $N>1$.
By the transience assumption and Lemma \ref{lemma}, $\gamma_x>0$, for all $x$.
 Thus, letting $n\to\infty$ in \eqref{key}, we obtain
\begin{equation}\label{finalupper-ii}
\frac1N\sum_{x=1}^N(2\omega(x)-1)\frac{1-\omega(x)\gamma_x}{1-\omega(x)}\ge1.
\end{equation}
Since \eqref{finalupper-ii} was derived under the assumption that $\gamma_x>0$, for all $x$, it follows
from \eqref{finalupper-ii} that $\frac1N\sum_{x=1}^N\frac{2\omega(x)-1}{1-\omega(x)}>1$ is a necessary condition for transience.
This is equivalent to $\frac1N\sum_{x=1}^N\frac{\omega(x)}{1-\omega(x)}>2$.

The analysis above has given necessary conditions for transience in parts (i) and (ii), and has also given the upper bound on $P_1(T_0=\infty)$ in part (i).
 We now turn to necessary conditions for recurrence in parts (i) and (ii), and to the lower bound on $P_1(T_0=\infty)$ in part (i).
For the time being, we make no assumption of transience or recurrence.
Under $P_1$, consider the random variable $D^x_{T_0\wedge T_n}$ appearing on the right hand side of \eqref{Doob}.
If $T_0<T_x$, then $D^x_{T_0\wedge T_n}=0$.
Let $\epsilon_{x,n}$ denote the probability  that after the first time the process
jumps rightward from $x$ to $x+1$, it does not  return to $x$ before reaching $n$. Similarly, for each $k>1$, let
$\epsilon^{(k)}_{x,n}$  denote the probability, conditioned on the process returning to $x$ from $x+1$
at least $k-1$ times before hitting $n$, that after the $k$-th time the process jumps rightward from $x$ to $x+1$,
it does not return to $x$ before reaching $n$.
Each time  the process jumps from  $x$ to $x+1$,
the current environment, call it $\omega''$, will satisfy $\omega''\le\omega'$, where $\omega'$ was the environment that was
in effect the previous time the process jumped from  $x$ to $x+1$.
Thus, by Proposition \ref{mono}, it follows that
$\epsilon^{(k)}_{x,n}\le\epsilon_{x,n}$, for $k>1$. In light of these
facts, it follows that the random variable $D^x_{T_0\wedge T_n}$  stochastically dominates
$I_xY_x$, where $I_x$ and $Y_x$ are independent random variables, with $I_x$ as in \eqref{IZ} and
with $Y_x$ satisfying
\begin{equation}
P(Y_x=k)=(1-\omega(x)+\omega(x)\epsilon_{x,n})(\omega(x)-\omega(x)\epsilon_{x,n})^{k-1},\  k\ge1.
\end{equation}
Thus,
\begin{equation}\label{DxTn}
\begin{aligned}
&E_1D^x_{T_n\wedge T_0}\ge EI_xY_x= \frac1{1-\omega(x)+\omega(x)\epsilon_{x,n}}P_1(T_x<T_0)\\
&\ge\frac1{1-\omega(x)+\omega(x)\epsilon_{x,n}}P_1(T_n<T_0), \ x=1,\cdots, n-1.
\end{aligned}
\end{equation}
From \eqref{Doob} and \eqref{DxTn}  it follows that for any $n$,
\begin{equation}\label{key2}
\frac1n\sum_{x=1}^{n-1}\frac{2\omega(x)-1}{1-\omega(x)+\omega(x)\epsilon_{x,n}}<1.
\end{equation}

Consider now part (i) of the theorem. In this case $\omega(x)=p$ for all $x$, and it follows that
$\epsilon_{x,n}$ depends only on $n-x$. It also follows that
$\lim_{(n-x)\to\infty}\epsilon_{x,n}=\gamma\equiv P_1(T_0=\infty)$.
If the process is recurrent, then $\gamma=0$.
Using these facts and letting $n\to\infty$ in   \eqref{key2}, we conclude
that $\frac{2p-1}{1-p}\le1$ is a necessary condition for recurrence. This is equivalent to $p\le\frac23$.
We also conclude that in the transient case,
$\frac{2p-1}{1-p+p\gamma}\le 1$ or equivalently,
\begin{equation}\label{uppergamma}
\gamma=P_1(T_0=\infty)\ge\frac{3p-2}p.
\end{equation}
This gives the lower bound on $P_1(T_0=\infty)$ in part (i).

Consider now part (ii) of the theorem and assume recurrence.
Then $\lim_{n\to\infty}\epsilon_{x,n}=0$.
Since $\omega(x)$ is periodic with period $N>1$,
the environment to the right of $x$ is the same as the environment to the right of $x+N$.
Thus,  $\epsilon_{x,n}=\epsilon_{x+N,n+N}$, and therefore   $\lim_{n-x\to\infty}\epsilon_{x,n}=0$.
Using these facts and letting $n\to\infty$ in \eqref{key2}, it follows that
a necessary condition for recurrence is
$\frac1N\sum_{x=1}^N\frac{2\omega(x)-1}{1-\omega(x)}\le 1$.
This is equivalent to $\frac1N\sum_{x=1}^N\frac{\omega(x)}{1-\omega(x)}\le 2$.

\hfill $\square$
\medskip

\noindent \bf Proof of Theorem 2.\rm\
By Lemma \ref{lemma}, the process is either recurrent for $\mathcal{P}^S$-almost every $\omega$ or transient for $\mathcal{P}^S$-almost
every $\omega$. Assume for now almost sure transience, and let $\omega$ be an environment  for which the process is transient.
For the fixed environment $\omega$, \eqref{key} holds; that is,
\begin{equation}\label{key-random}
P^\omega_1(T_n<T_0)\le\frac1n+\frac1n\sum_{x=1}^{n-1}(2\omega(x)-1)\frac{1-\omega(x)\gamma_x(\omega)}{1-\omega(x)}P^\omega_1(T_x<T_0),
\end{equation}
where $\gamma_x=\gamma_x(\omega)=P^\omega_{x+1}(T_x=\infty)$.
Note that $(2\omega(x)-1)\frac{1-\omega(x)\gamma_x(\omega)}{1-\omega(x)}$ can  be expressed in the form
$F(\theta^x\omega)$, where $\theta^x$ is the shift operator defined by $\theta^x(w)(z)=\omega(x+z)$.
For each $\omega$, we have $\lim_{z\to\infty}P_1^\omega(T_z<T_0)=\gamma_1(\omega)$.
Thus, letting $n\to\infty$ in \eqref{key-random} and using the stationarity and ergodicity, it follows
that
\begin{equation}\label{Exp1}
\mathcal{E}^S(2\omega(0)-1)\frac{1-\omega(0)\gamma_0(\omega)}{1-\omega(0)}\ge1.
\end{equation}
By the transience assumption and Lemma \ref{lemma},  $\gamma_0>0$ a.s. $\mathcal{P}^S$. Thus it follows from \eqref{Exp1}
that $\mathcal{E}^S\frac{2\omega(0)-1}{1-\omega(0)}>1$ is a necessary condition for transience.
This is equivalent to $\mathcal{E}^S\frac{\omega(0)}{1-\omega(0)}>2$.

Now assume almost sure recurrence. For any fixed environment $\omega$, \eqref{key2} holds; that is,
\begin{equation*}
\frac1n\sum_{x=1}^n\frac{2\omega(x)-1}{1-\omega(x)+\omega(x)\epsilon_{x,n}(\omega)}<1,
\end{equation*}
where $\epsilon_{x,n}(\omega)$ is the $P_1^\omega$-probability that after the first time the process jumps  rightward from $x$ to $x+1$, it does not return
to $x$ before reaching $n$.
Since $\epsilon_{x,n}(\omega)\le\epsilon_{x,x+M}(\omega)$, for $x\le n-M$, we have
\begin{equation}\label{key2-random}
\frac1n\sum_{x=1}^{n-M}\frac{2\omega(x)-1}{1-\omega(x)+\omega(x)\epsilon_{x,x+M}(\omega)}<1.
\end{equation}
By the stationarity and ergodicity,
\begin{equation}\label{conv}
\lim_{r\to\infty}\frac1r\sum_{x=1}^r\epsilon_{x,x+M}(\omega)=\mathcal{E}^S\epsilon_{0,M}\ \text{ a.s.},
\end{equation}
and by the recurrence assumption,
\begin{equation}\label{M}
\lim_{M\to\infty}\mathcal{E}^S\epsilon_{0,M}=0.
\end{equation}
From \eqref{conv} and \eqref{M} it follows that for each $\delta>0$ there exist $r_\delta$, $M_\delta$ such that
\begin{equation}\label{number}
\mathcal{P}^S(\frac1r\#\{x\in\{1,\cdots, r\}:\epsilon_{x,x+M}<\delta\}>1-\delta)>1-\delta, \ \text{for}\
r>r_\delta, M>M_\delta.
\end{equation}
By the stationarity and ergodicity again,
\begin{equation}\label{stat,erg}
\lim_{n\to\infty}\frac1n\sum_{x=1}^{n-M}\frac{2\omega(x)-1}{1-\omega(x)}=\mathcal{E}^S\frac{2\omega(0)-1}{1-\omega(0)}\ \text{a.s.}.
\end{equation}
From \eqref{key2-random}, \eqref{number} and \eqref{stat,erg}, we deduce  that
 a necessary condition for recurrence is
$\mathcal{E}^S\frac{2\omega(0)-1}{1-\omega(0)}\le1$. This is equivalent to
$\mathcal{E}^S\frac{\omega(0)}{1-\omega(0)}\le2$.
\hfill $\square$

\section{Proof of Theorem \ref{ball} }\label{3}
Since the random walk $\{X_n\}_{n=0}^\infty$ is recurrent if $p\in(\frac12,\frac23]$, we may assume that $p>\frac23$.
The proof begins  like the proof of the corresponding result in \cite{BS}.
Let
\begin{equation}
U_m^n=\#\{0\le k<T_n:X_k=m, X_{k+1}=m-1\}
\end{equation}
denote the number of times the process jumps from $m$ to $m-1$ before reaching $n$,
and let
\begin{equation}
K_n=\#\{0\le k\le T_n:X_k<0\}.
\end{equation}
An easy combinatorial argument yields
\begin{equation}
T_n=K_n-U_0^n+n+2\sum_{m=0}^n U_m^n.
\end{equation}
Since $p>\frac23$,     by Theorem \ref{t/r-det}  the
process is transient, and by Lemma \ref{lemma} one has
$\lim_{n\to\infty}X_n=\infty$ a.s.. Thus, $K_n$ and $U_0^n$ are
almost surely bounded, and we conclude that almost surely $T_n\sim
n+2\sum_{m=0}^nU_m^n$ as $n\to\infty$. We will show below that
there is a time-homogeneous  irreducible  Markov process $\{Y_0,
Y_1,\cdots\}$ such that for each $n$, the distribution of
$\{U_n^n,\cdots, U_0^n\}$ under $P_0$ is the same as the distribution of
$\{Y_0,  \cdots, Y_n\}$. By the transience of
$\{X_n\}_{n=0}^\infty$, $U_0^n$ converges to a limiting
distribution as $n\to\infty$. Thus, the same is true of $Y_n$,
which means that this Markov process is positive recurrent. Denote
its  stationary distribution on $Z^+$ by $\mu$. By the ergodic
theorem for irreducible Markov chains, we conclude that
\begin{equation}\label{ergodic}
\lim_{n\to\infty}\frac{T_n}n=1+2\sum_{k=0}^\infty k\mu(k)\le \infty.
\end{equation}

Now as is well-known \cite[p. 113]{Z}, for any $v\in[0,\infty)$, one has
$\lim_{n\to\infty}\frac{X_n}n=v$ a.s. if and only if $\lim_{n\to\infty}\frac{T_n}n=\frac1v$.
Thus in light of \eqref{ergodic}, the process $\{X_n\}_{n=0}^\infty$ is ballistic if and only if
$\sum_{k=0}^\infty k\mu(k)<\infty$, and the speed of the process is given by
$\frac1{1+2\sum_{k=0}^\infty k\mu(k)}$.

We now show
 that there is a time-homogeneous, irreducible  Markov process
$\{Y_0,Y_1,\cdots\}$ such that for each $n$, the distribution of
$\{U_n^n,\cdots, U_0^n\}$ under $P_0$  is the same as the distribution of
$\{Y_0,  \cdots, Y_n\}$, and we calculate its transition probabilities.
For this we need to define an IID sequence as follows.
Let $P$ denote probabilities for
an IID sequence $\{\zeta_n\}_{n=1}^\infty$  of geometric random variables satisfying
$P(\zeta_1=m)=\frac12^{m+1}$, $m=0,1,\cdots$.
Consider now the distribution under $P_0$  of the random variable $U_m^n$, conditioned on $U_{m+1}^n=j$.
 First consider the case $j=0$.
The process starts at 0 and eventually reaches $m$ for the first time. From $m$, it jumps right to $m+1$
 with probability $p$, and   jumps left to $m-1$  with
probability $1-p$.
 If the process jumps right to $m+1$, then it will never again find
itself at $m$. On the other hand, if it jumps left to $m-1$, then it will eventually return to $m$. When it returns
to $m$, it jumps left or right with probability $\frac12$. If it jumps right, then it will never again return to $m$, but
if it jumps left, then it will eventually return to $m$ again, where it again jumps left or right with probability $\frac12$.
This situation is repeated until the process finally jumps right to $m+1$, from which it will never return to $m$.
From the above description, it follows that $P_0(U_m^n=0|U_{m+1}^n=0)=p$ and $P_0(U_m^n=k|U_{m+1}^n=0)=(1-p)P(\zeta_1=k-1)$, for
$k\ge1$.

Now consider $j\ge1$. As above, the process starts at 0 and eventually reaches $m$ for the first time. From $m$, it jumps right to $m+1$
 with probability $p$, and   jumps left to $m-1$  with
probability $1-p$. If it jumps to the right, then since $j\ge1$ it will eventually return to $m$. Upon returning to $m$ it will
again jump right with probability $p$ and left with probability $1-p$. By the conditioning, it will jump back to $m$ from $m+1$ a total of
$j$ times. As long as the process
keeps jumping to the right when it is at $m$, the probability of jumping to the right the next time it is at $m$ will
still be $p$. But as soon as it jumps left from $m$, then whenever it again reaches $m$ it will jump right  with probability
$\frac12$. From this description, it follows that
$P_0(U_m^n=0|U_{m+1}^n=j)=p^{j+1}$ and
$P_0(U_m^n=k|U_{m+1}^n=j)=(1-p)\sum_{l=0}^jp^lP(\sum_{i=1}^{j+1-l}\zeta_i=k-1)$.

We have thus shown that  there exists such a time homogeneous Markov process
$\{Y_0,Y_1,\cdots\}$. Denote probabilities and expectations for this process starting from
$m\in\{0,1,\cdots\}$ by $\mathcal{P}_m$ and $\mathcal{E}_m$. We have also shown that the transition probabilities of this Markov
process are given by
\begin{equation}\label{transition}
p_{jk}\equiv \mathcal{P}(Y_1=k|Y_0=j)=
\begin{cases} p^{j+1},\  \text{if}\ k=0;\\
(1-p)\sum_{l=0}^jp^lP(\sum_{i=1}^{j+1-l}\zeta_i=k-1), \ \text{if}\ k\ge1.
\end{cases}
\end{equation}
As noted above, since $\{X_n\}_{n=0}^\infty$ is transient, the process $\{Y_n\}_{n=0}^\infty$ is positive recurrent. Furthermore, denoting its invariant measure
by $\mu$, the process is ballistic if  and only if
$\sum_{k=0}^\infty k\mu(k)<\infty$.

We now analyze $\sum_{k=0}^\infty k\mu(k)$.
Let $\sigma_m=\inf\{n\ge1: Y_n=m\}$.
By Khasminskii's representation \cite[chapter 5, section 4]{Du} of invariant probability measures, one has
\begin{equation}\label{Khas}
\mu_k=\frac{\mathcal{E}_1\sum_{m=1}^{\sigma_1}\chi_{k}(Y_m)}
{\mathcal{E}_1\sigma_1}.
\end{equation}
From \eqref{Khas}, it follows that
$\sum_{k=0}^\infty k\mu(k)<\infty$ if and only if
\begin{equation}\label{condition}
\mathcal{E}_1\sum_{m=1}^{\sigma_1}Y_m<\infty.
\end{equation}

 The generator of $\{Y_n\}_{n=0}^\infty$ is $P-I$, where
$Pu(j)=\sum_{k=0}^\infty p_{jk}u(k)$.
Thus, letting
\begin{equation}\label{martingale}
M^u_0=u(Y_0), \ \   M^u_n=u(Y_n)-\sum_{m=0}^{n-1} (Pu-u)(Y_m), \  n\ge1,
\end{equation}
it follows
that $\{M^u_n\}_{n=0}^\infty$  is a  martingale for all functions $u$ which satisfy
\begin{equation}\label{technical}
\sum_{k=0}^\infty p_{jk}|u(k)|<\infty,\ \text{ for all}\ j.
\end{equation}
Consider the function $u(k)=(k-1)^2$.
The ensuing calculation shows  that $u$ satisfies \eqref{technical}.
From \eqref{technical} and the fact that $E\zeta_1=1$ and $E\zeta_1^2=3$,
 it follows that
\begin{equation}\label{squared}
\begin{aligned}
&Pu(j)=p^{j+1}u(0)+(1-p)\sum_{l=0}^jp^lE(\sum_{i=1}^{j+1-l}\zeta_i)^2=\\
&p^{j+1}+
(1-p)\sum_{l=0}^jp^l\left((j+1-l)^2+2(j+1-l)\right).
\end{aligned}
\end{equation}
From \eqref{squared} and the fact that $(1-p)\sum_{l=0}^\infty lp^l=\frac p{1-p}$, it follows that
\begin{equation}\label{asympj2}
(Pu-u)(j)=(6-2\frac p{1-p})j+0(1),\ \text{as} \ j\to\infty.
\end{equation}
From Doob's optional sampling theorem,
\begin{equation}\label{optsamp}
\mathcal{E}_1M^u_{\sigma_1\wedge n}=\mathcal{E}_1M^u_0=0,\ n\ge0,
\end{equation}
and from \eqref{asympj2} and \eqref{martingale}, one has
\begin{equation}\label{final>}
\mathcal{E}_1M^u_{\sigma_1\wedge n}=\mathcal{E}_1(Y_{\sigma_1\wedge n}-1)^2-(6-2\frac p{1-p})\mathcal{E}_1\sum_{m=0}^{\sigma_1\wedge n-1}
Y_m+\ O(\mathcal{E}_1\sigma_1\wedge n)\
\text{as}\ n\to\infty.
\end{equation}
By the  positive recurrence of $\{Y_n\}_{n=0}^\infty$, which holds for all $p>\frac23$, one has $\mathcal{E}_1\sigma_1\wedge n\rightarrow \mathcal{E}_1\sigma_1<\infty$ as $n\to\infty$.
Thus \eqref{optsamp} and \eqref{final>} yield
\begin{equation}\label{mining}
\begin{aligned}
&\mathcal{E}_1(Y_{\sigma_1\wedge n}-1)^2-(6-2\frac p{1-p})\mathcal{E}_1\sum_{m=0}^{\sigma_1\wedge n-1}
Y_m=C_n,\ \text{for all}\ n,\\
&  \text{where}\ C_n=O(1),\ \text{as}\ n\to\infty.
\end{aligned}
\end{equation}
If $p>\frac34$, one has $6-2\frac p{1-p}<0$.
Letting
$n\to\infty$ in \eqref{mining}, we conclude
 that \eqref{condition} holds if $p>\frac34$, and thus the process $\{X_n\}_{n=0}^\infty$ is ballistic.

To calculate the speed of the process we will need to build on some calculations used in the proof of Proposition \ref{T0} below.
Thus we defer the proof of \eqref{explicit} until after the proof of Proposition \ref{T0}.

For the case $p\in(\frac23,\frac34)$, we will need the following proposition.
\begin{proposition}\label{T0}
Let $\{X_n\}_{n=0}^\infty$ be a random walk in a deterministic ``have your cookie and eat it'' environment
with $\omega(x)=p$, for all $x$.
Then
\begin{equation}
\begin{aligned}
&P_1(T_0=\infty)\ge\frac{3p-2}{2p-1}, \ \text{for}\ p\in(\frac23,\frac34);\\
&P_1(T_0=\infty)=\frac{3p-2}{2p-1}, \ \text{for}\ p\in[\frac34,1).
\end{aligned}
\end{equation}
\end{proposition}
\noindent \bf Remark.\rm\ We believe that equality holds in Proposition \ref{T0} for all $p\in(\frac23,1)$.

We will also need the following corollary, whose  proof will be pointed out during the proof
of Proposition \ref{T0}.
\begin{corollary}\label{T0cor}
Let $p\in(\frac23,1)$.
Let $\sigma_1=\inf\{m\ge1:Y_m=1\}$.
Then
$l\equiv\lim_{n\to\infty}E_1Y_{\sigma_1\wedge n}$ exists, and $l\ge1$, with equality if and only if
$P_1(T_0=\infty)=\frac{3p-2}{2p-1}$.
\end{corollary}

We now prove  non-ballisticity  in the case
 $p\in(\frac23,\frac34)$. Then we will prove Proposition \ref{T0}.
 We break the proof into two cases---$P_1(T_0=\infty)>\frac{3p-2}{2p-1}$ and $P_1(T_0=\infty)=\frac{3p-2}{2p-1}$.
The first case we can dispose of easily. (By the remark above, we expect that this case is actually vacuous.)
By Corollary \ref{T0cor}, we have $\lim_{n\to\infty}E_1Y_{\sigma_1\wedge n}>1$.
Using this along with the fact that $\lim_{n\to\infty}Y_{\sigma_1\wedge n}=Y_{\sigma_1}=1$ a.s.,  it follows that
\begin{equation}\label{squaredagain}
\lim_{n\to\infty}E_1(Y_{\sigma_1\wedge n}-1)^2=\infty.
\end{equation}
Now letting $n\to\infty$ in \eqref{mining}, and using \eqref{squaredagain} and the fact that $6-2\frac p{1-p}>0$, it follows
that \eqref{condition} does not hold.

We now treat the case $P_1(T_0=\infty)=\frac{3p-2}{2p-1}$. We take a completely different tack from the above.
In \cite{Z}, it was shown that for the random walks with cookies in that paper, one has
\begin{equation*}
\lim_{n\to\infty}\frac{X_n}n=\frac1u,\
 P_0-\text{a.s.},\ \text{where}\ u\equiv\sum_{j=1}^\infty P_0(T_{j+1}-T_j\ge j).
\end{equation*}
The reader can check that the same proof works for the  random walk in a ``have your cookie and eat it'' environment.
Thus, to prove nonballisticity in the case $p\in(\frac23,\frac34)$, it suffices to  show that $u=\infty$

For $m<j$, let $\nu_{m;j}=\inf\{k\ge T_j: X_k=m\}$. From the definition of the ``have your cookie and eat it'' environment,
on the event $\{\nu_{j-[\sqrt j~];j}<T_{j+1}\}$, after hitting $j-[\sqrt j]$ at time $\nu_{j-[\sqrt j];j}$, the process will eventually
find itself at $j-[\frac12\sqrt j]$, with the environment up to a distance $[\frac12\sqrt j]$ on either side of it being
a simple, symmetric random walk environment.
From $j-[\frac12\sqrt j]$, the process must exit the interval $[j-[\sqrt j], j]$ before time $T_{j+1}$.
By the scaling properties of the simple, symmetric random walk, there is
 a positive probability, independent of $j$, that it will take
at least $j$ steps to exit this interval, starting from its mid-point $j-[\frac12\sqrt j]$.
Thus, we conclude that for some $c>0$,
\begin{equation}
P_0(T_{j+1}-T_j\ge j)\ge c P_0(\nu_{j-[\sqrt j~];j}<T_{j+1}).
\end{equation}
Therefore to complete the proof it suffices to show that
\begin{equation}\label{only}
\sum_{j=1}^\infty P_0(\nu_{j-[\sqrt j~];j}<T_{j+1})=\infty.
\end{equation}

From the beginning of the proof of Theorem \ref{t/r-det}, recall that
$\sigma_x=\inf\{n\ge1: X_{n-1}=x, X_n=x-1\}$.
Define
\begin{equation}\label{rho}
\rho_x^{(j)}=\begin{cases} p, \ \text{if}\ \sigma_x>T_j,\\\frac12,\ \text{if}\ \sigma_x<T_j,\end{cases} \ \
\rho^{(\infty)}_x=\begin{cases} p, \ \text{if}\ \sigma_x=\infty,\\ \frac12,\ \text{if}\ \sigma_x<\infty.\end{cases}
\end{equation}
We introduce the notation  $P_0^{\{\rho_x^{(j)}\}}(\cdot)\equiv P_0(\cdot|\{\rho_x^{(j)}\}_{x=1}^{j-1})$.
We  can write
\begin{equation}\label{product}
\begin{aligned}
&P_0^{\{\rho_x^{(j)}\}}(\nu_{j-[\sqrt j];j}<T_{j+1})=\\
&P_0^{\{\rho_x^{(j)}\}}(\nu_{j-1;j}<T_{j+1})
\prod_{r=1}^{[\sqrt j]-1} P_0^{\{\rho_x^{(j)}\}}(\nu_{j-r-1;j}<T_{j+1}|\nu_{j-r;j}<T_{j+1}).
\end{aligned}
\end{equation}
Clearly,
\begin{equation}\label{clear}
P_0^{\{\rho_x^{(j)}\}}(\nu_{j-1;j}<T_{j+1})=1-p.
\end{equation}
The conditional probability   $P_0^{\{\rho_x^{(j)}\}}(\nu_{j-r-1;j}<T_{j+1}|\nu_{j-r;j}<T_{j+1})$
 is the probability that a simple random walk starting from
$j-r$ will hit $j-r-1$ before it hits $j+1$, where the probability of jumping to the right is
equal to $\frac12$ at the sites $x> j-r$ and is equal to $\rho^{(j)}_{j-r}$ at the site $j-r$.
This probability can be calculated explicitly using a standard difference equation with boundary conditions.
We leave it to the reader to check that one obtains
\begin{equation}\label{diffequ}
P_0^{\{\rho_x^{(j)}\}}(\nu_{j-r-1;j}<T_{j+1}|\nu_{j-r;j}<T_{j+1})=1-\frac{\rho_{j-r}^{(j)}}{(r+1)(1-\rho_{j-r}^{(j)})+\rho_{j-r}^{(j)}}.
\end{equation}
For every $\epsilon>0$, there exists an $x_\epsilon>0$ such that
$\log(1-x)>-(1+\epsilon)x$, for $0<x<x_\epsilon$.
Thus, from \eqref{rho}-\eqref{diffequ}, it follows that there exists an $r_\epsilon$ such that
\begin{equation}\label{bigequ}
\begin{aligned}
&P_0^{\{\rho_x^{(j)}\}}(\nu_{j-[\sqrt j];j}<T_{j+1})\ge\\
&(1-p)\left(\prod_{r=1}^{r_\epsilon}1-\frac{\rho_{j-r}^{(j)}}{r(1-\rho_{j-r}^{(j)})+\rho_{j-r}^{(j)}}\right)
\exp\left(-(1+\epsilon)\sum_{r=r_\epsilon+1}^{[\sqrt j]-1}\frac{\rho_{j-r}^{(j)}}{r(1-\rho_{j-r}^{(j)})+\rho_{j-r}^{(j)}}\right)\ge\\
&C_{\epsilon, p}\exp\left(-(1+\epsilon)\sum_{r=r_\epsilon+1}^{[\sqrt j]-1}\frac{\rho_{j-r}^{(j)}}{r(1-\rho_{j-r}^{(j)})}\right).
\end{aligned}
\end{equation}
Using \eqref{bigequ} and Jensen's inequality, we obtain
\begin{equation}\label{main}
\begin{aligned}
&P_0(\nu_{j-[\sqrt j];j}<T_{j+1})=E_0P_0^{\{\rho_x^{(j)}\}}(\nu_{j-[\sqrt j];j}<T_{j+1})\ge\\
& C_{\epsilon, p}E_0\exp\left(-(1+\epsilon)\sum_{r=r_\epsilon+1}^{[\sqrt j]-1}\frac{\rho_{j-r}^{(j)}}{r(1-\rho_{j-r}^{(j)})}\right)\ge\\
&C_{\epsilon, p}\exp\left(-(1+\epsilon)\sum_{r=r_\epsilon+1}^{[\sqrt j]-1}(E_0\frac{\rho_{j-r}^{(j)}}{1-\rho_{j-r}^{(j)}})\frac1r\right).
\end{aligned}
\end{equation}

Note that $\lim_{j\to\infty}\rho^{(j)}_x=\rho^{(\infty)}_x$. The distribution of $\rho^{(\infty)}_x$ is independent
of $x$ and is given by $P_0(\rho^{(\infty)}_x=p)=1-P_0(\rho^{(\infty)}_x=\frac12)=\gamma$, where
$\gamma\equiv P_1(T_0=\infty)$. Note that in fact the distribution of
$\rho^{(j)}_{j-r}$ depends only on $r$ and   converges weakly to the distribution
of $\rho_x^{(\infty)}$ as $r\to\infty$. We have $E_0\frac{\rho^{(\infty)}_x}{1-\rho^{(\infty)}_x}=\gamma\frac p{1-p}+1-\gamma$. In light of these
facts, and increasing $r_\epsilon$ if necessary, we have
\begin{equation}\label{delta}
\begin{aligned}
&\sum_{r=r_\epsilon+1}^{[\sqrt j]-1}(E_0\frac{\rho_{j-r}^{(j)}}{1-\rho_{j-r}^{(j)}})\frac1r
\le\\
&(\gamma\frac p{1-p}+1-\gamma+\epsilon)\log\sqrt j.
\end{aligned}
\end{equation}
From  \eqref{main} and \eqref{delta}, we obtain
\begin{equation}\label{final}
P_0(\nu_{j-[\sqrt j];j}<T_{j+1})\ge \frac{C_{\epsilon, p}}{j^{\frac12(1+\epsilon)(\gamma\frac p{1-p}+1-\gamma+\epsilon)}}.
\end{equation}
Since $\epsilon$ can be chosen arbitrarily small, it follows from \eqref{final}
that \eqref{only} holds if $\frac12(\gamma\frac p{1-p}+1-\gamma)<1$, or equivalently, if
$\gamma<\frac{1-p}{2p-1}$.
But we have assumed that  $\gamma= \frac{3p-2}{2p-1}$.
Thus, \eqref{only} holds if $\frac{3p-2}{2p-1}<\frac{1-p}{2p-1}$, which is equivalent to $p<\frac34$.
This completes the proof of non-ballisticity in the case $p\in(\frac23,\frac34)$.

We now prove Proposition \ref{T0}, after which  we will prove \eqref{explicit}, which as already noted requires some of the calculations in the proof of
the proposition.

\bigskip

\bf \noindent Proof of Proposition \ref{T0}.\rm\
Since for each $n$, $\{U_n^n,\cdots
U_1^n\}$ under $P_0$ is a Markov process with transition
probabilities  $p_{jk}$ as above in \eqref{transition} and with
stationary distribution $\mu$, it follows that
$\mu_0=\lim_{n\to\infty}P_0(U_1^n=0)$. On the other hand
$P_1(T_0=\infty)=\lim_{n\to\infty}P_0(U_1^n=0)$. Thus,
$P_1(T_0=\infty)=\mu_0$.
By the invariance of $\mu$ and \eqref{transition}, we have $\mu_0=\sum_{j=0}^\infty p_{j0}\mu_j=\sum_{j=0}^\infty p^{j+1}\mu_j$. Thus,
\begin{equation}\label{mu0}
P_1(T_0=\infty)=\sum_{j=0}^\infty p^{j+1}\mu_j.
\end{equation}
To prove the proposition, we will evaluate the right hand side of \eqref{mu0}.

Consider the martingale $\{M^u_n\}_{n=0}^\infty$ defined in \eqref{martingale}, with $u(k)=k-1$.
Using the fact that $E\zeta_1=1$, one has from \eqref{transition}
\begin{equation}
\begin{aligned}
&Pu(j)=\sum_{k=0}^\infty p_{jk}u(k)=
-p^{j+1}+(1-p)\sum_{l=0}^jp^lE(\sum_{i=1}^{j+1-l}\zeta_i)=\\
&-p^{j+1}+(1-p)\sum_{l=0}^jp^l(j+1-l)=\\
&-p^{j+1}+(j+1)(1-p^{j+1})-\frac p{1-p}(1-(j+1)p^j+jp^{j+1}).
\end{aligned}
\end{equation}
After regrouping terms, one obtains
\begin{equation}\label{arith}
(Pu-u)(j)=\frac1{1-p}\left(2-3p+p^{j+1}(2p-1)\right).
\end{equation}
By Doob's optional stopping theorem, $\{M_n^u\}_{n=0}^\infty$ satisfies \eqref{optsamp}, where $\sigma_1=\inf\{m\ge1:Y_m=1\}$.
Using this along with \eqref{arith} and the definition of $\{M^u_n\}_{n=0}^\infty$, we obtain
\begin{equation}\label{equal}
\mathcal{E}_1Y_{\sigma_1\wedge n}=1+ \mathcal{E}_1\sum_{m=0}^{\sigma_1\wedge n-1}
\frac1{1-p}\left(2-3p+p^{Y_m+1}(2p-1)\right).
\end{equation}
By Khasminskii's  representation of invariant probability measures \cite[chapter 5, section 4]{Du},
one has
\begin{equation}\label{khasagain}
\begin{aligned}
&\lim_{n\to\infty}\mathcal{E}_1\sum_{m=0}^{\sigma_1\wedge n-1}
\frac1{1-p}\left(2-3p+p^{Y_m+1}(2p-1)\right)=\\
&(\mathcal{E}_1\sigma_1)
\sum_{j=0}^\infty\frac1{1-p}\left(2-3p+p^{j+1}(2p-1)\right)\mu_j.
\end{aligned}
\end{equation}
This shows that
 the right hand side of \eqref{equal} converges as $n\to\infty$; thus,  so does the left hand side. Since
$\sigma_1<\infty$ a.s., and since $Y_{\sigma_1\wedge n}$ is bounded from below and $Y_{\sigma_1}=1$, it follows that the
left hand side of \eqref{equal} converges
to a  limit that is greater than or equal to 1 as $n\to\infty$.
Thus, the right hand side of \eqref{khasagain} is nonnegative; equivalently,
$\sum_{j=0}^\infty p^{j+1}\mu_j\ge\frac{3p-2}{2p-1}$. From this and \eqref{mu0}, it follows that
$P_1(T_0=\infty)\ge\frac{3p-2}{2p-1}$. The above calculation also shows that
Corollary \ref{T0cor} holds.
For the proof of non-ballisticity given above in the case $p\in(\frac23,\frac34)$, this is all we need. We now continue and show that
$P_1(T_0=\infty)=\frac{3p-2}{2p-1}$, if $p\ge\frac34$. This will be needed for the proof of \eqref{explicit}.

First consider the case $p=\frac34$. If on the contrary, one had
$P_1(T_0=\infty)>\frac{3p-2}{2p-1}$, then it would follow from Corollary \ref{T0cor} that
$\lim_{n\to\infty}E_1Y_{\sigma_1\wedge n}>1$, and as noted in the proof of Theorem \ref{ball}
(see \eqref{squaredagain}), this would
imply that $\lim_{n\to\infty}E_1(Y_{\sigma_1\wedge n}-1)^2=\infty$.
Using this and letting $n\to\infty$ in \eqref{mining} would lead to a contradiction
since $6-2\frac p{1-p}=0$.

Now consider the case  $p>\frac34$.
Let $u(k)=\exp(-\gamma (k-1))$, where $\gamma>0$. Since $\mu$ is the invariant probability measure, we have
\begin{equation}\label{invariant}
\sum_{j=0}^\infty Pu(j)\mu_j=\sum_{j=0}^\infty
u(j)\mu_j=\sum_{j=0}^\infty\exp(-\gamma(j-1))\mu_j.
\end{equation}
Using the fact that $E\exp(-\gamma\zeta_1)=(2-\exp(-\gamma))^{-1}$, we have
from \eqref{transition}
\begin{equation}\label{long}
\begin{aligned}
&Pu(j)=\sum_{k=0}^\infty p_{jk}\exp(-\gamma (k-1))=p^{j+1}\exp(\gamma)+\\
&\sum_{k=1}^\infty(1-p)\sum_{l=0}^jp^lP(\sum_{i=1}^{j+1-l}\zeta_i=k-1)\exp(-\gamma(k-1))=\\
&p^{j+1}\exp(\gamma)+(1-p)\sum_{l=0}^jp^lE\exp(-\gamma(\sum_{i=1}^{j+1-l}\zeta_i))=\\
&p^{j+1}\exp(\gamma)+(1-p)\sum_{l=0}^jp^l(\frac1{2-\exp(-\gamma)})^{j+1-l}=\\
&p^{j+1}\exp(\gamma)+(\frac1{2-\exp(-\gamma)})^{j+1}(1-p)\frac{1-(p(2-\exp(-\gamma)))^{j+1}}{1-p(2-\exp(-\gamma))}=\\
&p^{j+1}\exp(\gamma)+\frac{1-p}{p(2-\exp(-\gamma))-1}(p^{j+1}-(\frac1{2-\exp(-\gamma)})^{j+1}).
\end{aligned}
\end{equation}
Substituting for $Pu$ from the right hand side of \eqref{long} into the left hand side of
\eqref{invariant},
and then multiplying both sides of the resulting equation by $p(2-\exp(-\gamma))-1$,
we obtain
\begin{equation}\label{basic}
\begin{aligned}
&(2p-1)(\exp(\gamma)-1)\sum_{j=0}^\infty p^{j+1}\mu_j=(1-p)\sum_{j=0}^\infty(\frac1{2-\exp(-\gamma)})^{j+1}\mu_j+\\
&\left((2p-1)\exp(\gamma)-p\right)\sum_{j=0}^\infty\exp(-\gamma j)\mu_j.
\end{aligned}
\end{equation}
 Differentiating \eqref{basic} with respect to $\gamma$
 gives
\begin{equation}\label{longlong}
\begin{aligned}
&(2p-1)\exp(\gamma)\sum_{j=0}^\infty p^{j+1}\mu_j=
\left(p-(2p-1)\exp(\gamma)\right)\sum_{j=0}^\infty\exp(-\gamma j)j\mu_j-\\
&(1-p)\exp(-\gamma)(2-\exp(-\gamma))^{-2}
\sum_{j=0}^\infty (\frac1{2-\exp(-\gamma)})^jj\mu_j+\\
&(2p-1)\exp(\gamma)\sum_{j=0}^\infty\exp(-\gamma j)\mu_j-\\
&(1-p)\exp(-\gamma)(2-\exp(-\gamma))^{-2}\sum_{j=0}^\infty(\frac1{2-\exp(-\gamma)})^j\mu_j=\\
&I-II+III-IV.
\end{aligned}
\end{equation}
Letting $\gamma\to0$, the left hand side of \eqref{longlong} converges to $(2p-1)\sum_{j=0}^\infty p^{j+1}\mu_j$ and
  III$-$IV converges to $3p-2$.
Since $p>\frac34$, it follows from  Theorem \ref{ball} that   $\sum_{j=0}^\infty j\mu_j<\infty$. Thus,
I$-$II converges to $0$ as $\gamma\to0$, and consequently  one obtains $\sum_{j=0}^\infty p^{j+1}\mu_j=\frac{3p-2}{2p-1}$.
From this and \eqref{mu0} it follows that $P_1(T_0=\infty)=\frac{3p-2}{2p-1}$.
This concludes the proof of Proposition \ref{T0}.
(If one could show that I$-$II is non-positive for small $\gamma>0$, then one could  conclude from the above
calculation that
$P_1(T_0=\infty)\le\frac{3p-2}{2p-1}$, which in conjunction with the reverse inequality  proved above, would prove
that $P_1(T_0=\infty)=\frac{3p-2}{2p-1}$.
However our analysis of I$-$II proved inconclusive.)

\hfill $\square$

We now turn to the proof of \eqref{explicit}. As was noted in the paragraph following \eqref{ergodic}, the speed of the process
is given by $\frac1{1+2\sum_{k=0}^\infty k\mu_k}$.
Thus, to prove \eqref{explicit}, we need to show that $\sum_{k=0}^\infty k\mu_k\le\frac{2-2p}{4p-3}$, with equality if
$p\in[\frac{10}{11},1)$.
The calculations are in the spirit of several of the previous calculations, except that they are more tedious. We skip some of the intermediate
steps. First we show that
\begin{equation}\label{secondmoment}
\sum_{k=0}^\infty k^2\mu_k<\infty, \ \text{if}\ p\in(\frac{10}{11},1).
\end{equation}
Let $u(k)=(k-1)^3$. From \eqref{transition} we have
\begin{equation}\label{genagain}
Pu(j)=-p^{j+1}+(1-p)\sum_{l=0}^jp^lE(\sum_{i=1}^{j+1-l}\zeta_i)^3.
\end{equation}
We have already noted that $E\zeta_1=1$ and $E\zeta_1^2=3$. One can check that $E\zeta_1^3=13$. Also, one has
$(1-p)\sum_{l=0}^\infty lp^l=\frac p{1-p}$ and $(1-p)\sum_{l=0}^\infty l^2p^l=\frac p{1-p}+2(\frac p{1-p})^2$.
Using these facts, one calculates that
\begin{equation}\label{asympj3}
(Pu-u)(j)=j^2\left(30-\frac {3p}{1-p}\right)+j\left(34-49\frac p{1-p}+2(\frac p{1-p})^2\right)+O(1), \ \text{as}\ j\to\infty.
\end{equation}
Similar to \eqref{condition}, one has
from \eqref{Khas} that
$\sum_{k=0}^\infty k^2\mu(k)<\infty$ if and only if
\begin{equation}\label{conditionagain}
\mathcal{E}_1\sum_{m=1}^{\sigma_1}Y^2_m<\infty.
\end{equation}
Applying the argument from \eqref{martingale} to \eqref{mining}, but using the function $u(k)=(k-1)^3$ instead of the function $u(k)=(k-1)^2$,  using
\eqref{asympj3} instead of \eqref{asympj2}, and noting that $30-\frac {3p}{1-p}<0$ if $p>\frac{10}{11}$,
 one concludes that \eqref{conditionagain} holds if $p>\frac{10}{11}$. This proves \eqref{secondmoment}

We now prove the equality in \eqref{explicit} for $p\in[\frac{10}{11},1)$. The case $p=\frac{10}{11}$ uses the same type of argument
as that used to prove that $P_1(T_0=\infty)=\frac{3p-2}{2p-1}$ in the case $p=\frac34$. One uses the function $u(k)=(k-1)^3$ and the
calculations in the previous paragraph rather than the function $u(k)=(k-1)^2$. We leave the details to the reader. For the case
$p>\frac{10}{11}$, we
differentiate \eqref{longlong} in $\gamma$  and set $\gamma=0$. Two of the terms that result in this tedious
calculation  are $\pm(1-p)\sum_{k=0}^\infty k^2\mu_k$. By \eqref{secondmoment},
these terms pose
no problem. One obtains
\begin{equation}\label{longcal}
(2p-1)\sum_{j=0}^\infty p^{j+1}\mu_j=2-p+(-8p+6)\sum_{j=0}^\infty j\mu_j.
\end{equation}
In the proof of  Proposition \ref{T0}, it was shown that  $\sum_{j=0}^\infty p^{j+1}\mu_j=\frac{3p-2}{2p-1}$. Thus we conclude from \eqref{longcal} that
$\sum_{j=0}^\infty j\mu_j=\frac{2-2p}{4p-3}$.

We now turn to the case $p\in(\frac34,\frac{10}{11})$.
The argument from \eqref{arith} until the end of that paragraph, applied to $u(k)=(k-1)^2$ instead of to $u(k)=k-1$, shows that
\begin{equation}\label{khasineq}
\sum_{j=0}^\infty (Pu-u)(j)\mu_j\ge0.
\end{equation}
A very tedious calculation starting from the formula for $Pu(j)$ on the second line of \eqref{squared} reveals that
\begin{equation}\label{tedious}
(Pu-u)(j)=j(6-\frac {2p}{1-p})+\frac1{1-p}(-5p+\frac{2p^2}{1-p}-\frac{2p^{j+3}}{1-p}+2+p^{j+1}).
\end{equation}
Substituting \eqref{tedious} into \eqref{khasineq} and again using the fact that $\sum_{j=0}^\infty p^{j+1}\mu_j=\frac{3p-2}{2p-1}$,
one concludes after a bit of algebra that $\sum_{j=0}^\infty j\mu_j\le\frac{2-2p}{4p-3}$.
\hfill $\square$

\medskip

\bf \noindent Acknowledgment.\rm\ The author thanks a referee for the careful reading given to the paper and for  comments
that led to its improvement.

\end{document}